# NEW FRACTIONAL DERIVATIVES WITH NONLOCAL AND NON-SINGULAR KERNEL: THEORY AND APPLICATION TO HEAT TRANSFER MODEL

[1]Abdon Atangana and [2,3]Dumitru Baleanu

Institute for Groundwater Studies, Faculty of Natural and Agricultural Sciences,
University of the Free State, 9301, Bloemfontein, South Africa
Email: abdonatangana@yahoo.fr
[2]Cankaya University, Department of Mathematics and Computer Sciences, Ankara,
Turkey, email: dumitru@cankaya.edu.tr
[3]Institute of Space Sciences, Magurele-Bucharest, Romania

*In this manuscript we proposed a new fractional derivative with non-local and no-singular kernel. We presented some useful properties of the new derivative and applied it to solve the fractional heat transfer model.*

## 1. Introduction

Recently, a new derivative was launched by Caputo and Fabrizio [1] and it was followed by some related theoretical and applied results (see for example Refs.[2-4] and the references therein). We recall that the existing fractional derivatives have been used in many real world problems with great success (see for example Refs. [7-14] and the references therein) but still there are many thinks to be done in this direction.

**Definition 1:** [1] *Let* $f \in H^1(a,b), b > a, \alpha \in [0,1]$ *then, the definition of the new Caputo fractional derivative is:*

$$D_t^\alpha \big(f(t)\big) = \frac{M(\alpha)}{1-\alpha} \int_b^t f'(x) \exp\left[-\alpha \frac{t-x}{1-\alpha}\right] dx, \tag{1}$$

where $M(\alpha)$ denotes a normalization function obeying $M(0) = M(1) = 1$. However, if the function does not belongs to $H^1(a,b)$ then, the derivative has the form

$$D_t^\alpha \big(f(t)\big) = \frac{\alpha M(\alpha)}{1-\alpha} \int_b^t (f(t) - f(x)) \exp\left[-\alpha \frac{t-x}{1-\alpha}\right] dx. \tag{2}$$

If $\sigma = \frac{1-\alpha}{\alpha} \in [0,\infty), \alpha = \frac{1}{1+\sigma} \in [0,1]$, then the above equation assumes the form

$$D_t^\sigma(f(t)) = \frac{N(\sigma)}{\sigma} \int_b^t f'(x) \exp\left[-\frac{t-x}{\sigma}\right] dx, \qquad N(0) = N(\infty) = 1.$$

The aim of [1] was to introduce of a new derivative with exponential kernel. Its anti-derivative was reported in [2] and it was found to be the average of a given function. The derivative introduced in [1] cannot produce the original function when alpha is one. However this issue was, so far, independently solved in [5] and [6], respectively. We believe that the main message presented in [1] was to find a way to describe even better the dynamics of systems with memory effect. For a given data we ask the following question: what is the most accurate kernel which better describe it? We suggest a possible answer in the following sections.

## 2. New derivatives with nonlocal kernel

We recall that the Mittag-Leffler function is the solution of the following fractional ordinary differential equation [14-16]:

$$\frac{d^\alpha y}{dx^\alpha} = ay, 0 < \alpha < 1. \tag{3}$$

The Mittag-Leffler function and its generalized versions are therefore considered as nonlocal functions. Let us consider the following generalized Mittag-Leffler function:

$$E_\alpha(-t^\alpha) = \sum_{k=0}^{\infty} \frac{(-t)^{\alpha k}}{\Gamma(\alpha k + 1)} \tag{4}$$

The Taylor series of $\exp(-(t-y))$ at the point t is given by:

$$\exp(-a(t-y)) = \sum_{k=0}^{\infty} \frac{(-a(t-y))^k}{k!}. \tag{5}$$

If we chose $a = \frac{\alpha}{1-\alpha}$ and replace the above expression into Caputo-Fabrizio derivative we conclude that

$$D_t^\alpha(f(t)) = \frac{M(\alpha)}{1-\alpha} \sum_{k=0}^{\infty} \frac{(-a)^k}{k!} \int_b^t \frac{df(y)}{dy} ((t-y))^k dy. \tag{6}$$

To solve the problem of non-locality, we derive the following expression.

In equation (6), we replace $k!$ by $\Gamma(\alpha k + 1)$ also $(t-y)^k$ is replaced by $(t-y)^{\alpha k}$ to obtain:

$$D_t^\alpha(f(t)) = \frac{M(\alpha)}{1-\alpha} \sum_{k=0}^{\infty} \frac{(-a)^k}{\Gamma(\alpha k + 1)} \int_b^t \frac{df(y)}{dy} ((t-y))^{\alpha k} dy.$$

Thus, the following derivative is proposed.

**Definition 2:** *Let $f \in H^1(a,b), b > a, \alpha \in [0,1]$ then, the definition of the new fractional derivative is given as:*

$$^{ABC}_{\ \ b}D^\alpha_t(f(t)) = \frac{B(\alpha)}{1-\alpha}\int_b^t f'(x) E_\alpha\left[-\alpha\frac{(t-x)^\alpha}{1-\alpha}\right]dx. \tag{7}$$

Of course B(α) has the same properties as in Caputo and Fabrizio case. The above definition will be helpful to discuss real world problems and it also will have a great advantage when using the Laplace transform to solve some physical problem with initial condition. However, when alpha is 0 we do not recover the original function except when at the origin the function vanishes. To avoid this issue, we propose the following definition.

**Definition 3:** *Let $f \in H^1(a,b), b > a, \alpha \in [0,1]$ then, the definition of the new fractional derivative is given as:*

$$^{ABR}_{\ \ b}D^\alpha_t(f(t)) = \frac{B(\alpha)}{1-\alpha}\frac{d}{dt}\int_b^t f(x) E_\alpha\left[-\alpha\frac{(t-x)^\alpha}{1-\alpha}\right]dx. \tag{8}$$

Equations (7) and (8) have a non-local kernel. Also in equation (7) when the function is constant we get zero.

### 2.1. Properties of the new derivatives

In this section, we start by presenting the relation between both derivatives with Laplace transform. By simple calculation we conclude that

$$\mathcal{L}\{^{ABR}_{\ \ 0}D^\alpha_t(f(t))\}(p) = \frac{B(\alpha)}{1-\alpha}\frac{p^\alpha \mathcal{L}\{f(t)\}(p)}{p^\alpha + \frac{\alpha}{1-\alpha}} \tag{9}$$

and

$$\mathcal{L}\{^{ABC}_{\ \ 0}D^\alpha_t(f(t))\}(p) = \frac{B(\alpha)}{1-\alpha}\frac{p^\alpha \mathcal{L}\{f(t)\}(p) - p^{\alpha-1}f(0)}{p^\alpha + \frac{\alpha}{1-\alpha}}, \tag{10}$$

respectively.

The following theorem can therefore be established.

**Theorem 1:** *Let $f \in H^1(a,b), b > a, \alpha \in [0,1]$ then, the following relation is obtained*

$$^{ABC}_{\ \ 0}D^\alpha_t(f(t)) = {^{ABR}_{\ \ 0}}D^\alpha_t(f(t)) + H(t) \tag{11}$$

**Proof:** By using the definition (11) and the Laplace transform applied on both sides we obtain easily the following result:

$$\mathcal{L}\{^{ABC}_{\ \ 0}D^\alpha_t(f(t))\}(p) = \frac{B(\alpha)}{1-\alpha}\frac{p^\alpha \mathcal{L}\{f(t)\}(p)}{p^\alpha + \frac{\alpha}{1-\alpha}} - \frac{p^{\alpha-1}f(0)}{p^\alpha + \frac{\alpha}{1-\alpha}}\frac{B(\alpha)}{1-\alpha}. \tag{12}$$

Following equation (9) we have:

$$\mathcal{L}\{^{ABC}_{\ \ 0}D^\alpha_t(f(t))\}(p) = \mathcal{L}\{^{ABR}_{\ \ 0}D^\alpha_t(f(t))\}(p) - \frac{p^{\alpha-1}f(0)}{p^\alpha + \frac{\alpha}{1-\alpha}}\frac{B(\alpha)}{1-\alpha}. \tag{13}$$

Applying the inverse Laplace on both sides of equation (13) we obtain:

$$^{ABC}_0D^\alpha_t(f(t)) = {}^{ABR}_0D^\alpha_t(f(t)) - \frac{B(\alpha)}{1-\alpha}f(0)E_\alpha\left(-\frac{\alpha}{1-\alpha}t^\alpha\right). \tag{14}$$

This completes the proof.

**Theorem 2**: *Let f be a continuous function on a closed interval [a, b]. Then the following inequality is obtained on [a, b]*

$$\left\|{}^{ABR}_0D^\alpha_t(f(t))\right\| < \frac{B(\alpha)}{1-\alpha}K, \quad \|h(t)\| = \max_{a\leq t\leq b}|h(t)|. \tag{15}$$

**Proof:**

$$\left\|{}^{ABR}_0D^\alpha_t(f(t))\right\| = \left\|\frac{B(\alpha)}{1-\alpha}\frac{d}{dt}\int_0^t f(x)E_\alpha\left[-\alpha\frac{(t-x)^\alpha}{1-\alpha}\right]dx\right\| < \frac{B(\alpha)}{1-\alpha}\left\|\frac{d}{dt}\int_0^t f(x)dx\right\| = \frac{B(\alpha)}{1-\alpha}\|f(x)\|.$$

Then taking K to be $\|f(x)\|$ the proof is completed.

**Theorem 3**: *The A.B. derivative in Riemann and Caputo sense possess the Lipschitz condition, that is to say, for a given couple function f and h, the following inequalities can be established:*

$$\left\|{}^{ABR}_0D^\alpha_t(f(t)) - {}^{ABR}_0D^\alpha_t(h(t))\right\| \leq H\|f(t) - h(t)\| \tag{16}$$

and also

$$\left\|{}^{ABC}_0D^\alpha_t(f(t)) - {}^{ABC}_0D^\alpha_t(h(t))\right\| \leq H\|f(t) - h(t)\|. \tag{17}$$

We present the proof of (16) as the proof of (17) can be obtained similarly.

**Proof:**

$$\left\|{}^{ABR}_0D^\alpha_t(f(t)) - {}^{ABR}_0D^\alpha_t(h(t))\right\| =$$

$$\left\|\frac{B(\alpha)}{1-\alpha}\frac{d}{dt}\int_0^t f(x)E_\alpha\left[-\alpha\frac{(t-x)^\alpha}{1-\alpha}\right]dx - \frac{B(\alpha)}{1-\alpha}\frac{d}{dt}\int_0^t h(x)E_\alpha\left[-\alpha\frac{(t-x)^\alpha}{1-\alpha}\right]dx\right\|.$$

Using the Lipschitz condition of the first order derivative, we can find a small positive constant such that:

$$\left\|{}^{ABR}_0D^\alpha_t(f(t)) - {}^{ABR}_0D^\alpha_t(f(t))\right\| < \frac{B(\alpha)\theta_1}{1-\alpha}E_\alpha\left[-\alpha\frac{t^\alpha}{1-\alpha}\right]\left\|\int_0^t f(x)dx - \int_0^t h(x)dx\right\| \tag{18}$$

and then the following result is obtained:

$$\left\|{}^{ABR}_0D^\alpha_t(f(t)) - {}^{ABR}_0D^\alpha_t(f(t))\right\| < \frac{B(\alpha)\theta_1}{1-\alpha}E_\alpha\left[-\alpha\frac{t^\alpha}{1-\alpha}\right]\|f(x) - h(x)\|t$$

$$= H\|f(x) - h(x)\|,$$

which produces the requested result.

*Let $f$ be an n-times differentiable with natural number and $f^{(k)}(0) = 0, k = 1,2,3, \ldots \ldots n$, then by inspection we obtain*

$$^{ABC}_{\ \ 0}D^{\alpha}_t \left(\frac{d^n f(t)}{dt^n}\right) = \frac{d^n}{dt^n}\left(^{ABR}_{\ \ 0}D^{\alpha}_t(f(t))\right). \tag{19}$$

Now, we can easily prove by taking the inverse Laplace transform and using the convolution theorem that *the following time fractional ordinary differential equation:*

$$^{ABC}_{\ \ 0}D^{\alpha}_t(f(t)) = u(t) \tag{20}$$

*has a unique solution, namely*

$$f(t) = \frac{1-\alpha}{B(\alpha)}u(t) + \frac{\alpha}{B(\alpha)\Gamma(\alpha)}\int_0^t u(y)(t-y)^{\alpha-1}dy.$$

**Definition 4**: The fractional integral associate to the new fractional derivative with non-local kernel is defined as:

$$^{AB}_{\ a}I^{\alpha}_t\{f(t)\} = \frac{1-\alpha}{B(\alpha)}f(t) + \frac{\alpha}{B(\alpha)\Gamma(\alpha)}\int_a^t f(y)(t-y)^{\alpha-1}dy. \tag{21}$$

When α is zero we recover the initial function and if also α is 1, we obtain the ordinary integral.

### 3. Application to thermal science: A new heat transfer model

The analyses of the experimental data coming from the investigation of dynamics of complex systems is still an interesting open problem. Therefore, new methods and techniques are still to be discovered to and apply with more success to gave an even better description of the dynamics of real world problems. Particularly, finding new derivatives suggested in this manuscript is because of the necessity of employing a better model portraying the behaviour of orthodox viscoelastic materials, thermal medium and other. The new approach is able to portray material heterogeneities and some structure or media with different scales. The non-locality of the new kernel allows better description of the memory within structure and media with different scales. In addition of this, we also rely that this new derivative can play a specific role in the study of macroscopic behavior of some materials, related with nonlocal exchanges, which are predominant in defining the properties of the material [1].Thus, this new derivative will be very useful in describing many complex problems in thermal science. We recall here that another derivative was introduced in [17-18] with the aim of solving some problems within the scope of thermal science.

To describe the time rate of heat transfer through a material with different scale or heterogeneous, we propose a new law of heat conduction which will be refereed as fractional Fourier's law. The fractional Fourier law states that time rate of heat transfer via material with different scale is proportional to negative gradient in temperature and area at right angles to that gradient via which the heat flows and is given as:

$$^{ABC}_{0}D_t^\alpha Q = -k \oint_S \nabla T \cdot \overrightarrow{dA}. \tag{22}$$

Here $^{ABC}_{0}D_t^\alpha Q(r,t)$ is the amount of heat transferred within material with different scale per unit time. Now, for heterogeneous material of one dimensional geometry between two endpoints at constant temperature, produced a new heat flow rate as:

$$\frac{dQ}{dt} = -kA \, ^{ABR}_{0}D_r^\alpha T. \tag{23}$$

For instance in a cylindrical heterogeneous shells like a pipes, the heat conduction via an heterogeneous shell will be determined from internal radius $r_1$ and the external radius $r_2$, the length, $L$, and the difference between inner and outer wall we have the following

$$\frac{dQ}{dt} = 2\pi l r k A. \tag{24}$$

Rearranging and apply the Laplace transform on both sides, then applying the inverse Laplace transform, we obtain:

$$T_1 - T_2 = \frac{\dot{Q}}{2\pi lk}\left\{1 + \frac{\alpha}{1-\alpha}r_2^{-1+\alpha}\left\{\alpha\text{HypergeometricPFQ}\left[\{1,1,1-\alpha\},\{2,2\},\frac{r_1}{r_2}\right]r_1 - \left\{\text{HarmonicNumber}[\alpha] + \ln\left[\frac{r_1}{r_2}\right]\right\}r_2\right\}\right\}. \tag{25}$$

**Conclusions**

The aim of this manuscript was to suggest new derivatives with non-local and non-singular kernel. To achieve this goal, we make use the generalized Mittag-Leffler function to build the nonlocal kernel. One derivative is based upon the Caputo viewpoint and the second on the Riemann-Liouville approach. We derive the fractional integral associate using the Laplace transform operator. The new derivative was used to model the flow of heat in material with different scale and also those with heterogeneous media.

Nomenclature

T- temperature (C)

L- the length (m)

K- is a constant

R- is the radius (m)